\documentclass[11pt,openbib]{article}
\usepackage{latexsym}
\usepackage{amssymb,amsbsy,amsmath,amsfonts,amssymb,amscd}
\usepackage{graphics,pstricks}
\newcommand{\mZ}{{\mathbb Z}}               
\newcommand{\mR}{{\mathbb R}}               
\newcommand{\mC}{{\mathbb C}}               
\newcommand{\mT}{{\mathbb T}}               

\newcommand{\R}{\mathbb{R}}

\newcommand{\dee}{\mathop{\! \, \mathrm{d} \!}\nolimits}
\newcommand{\comp}{\raisebox{0pt}{$\scriptstyle\circ \, $}}

\newcommand{\onehalf}{\mbox{$\frac{\scriptstyle 1}{\scriptstyle 2}\,$}}

\newcommand{\lefthook}{\mbox{$\, \rule{8pt}{.5pt}\rule{.5pt}{6pt}\, \, $}}

\begin{document}
\begin{center}
{\Large \textbf{Geometry of KAM tori for nearly}}
\mbox{}\vspace{.05in} \\
{\Large \textbf{integrable Hamiltonian systems}}
\mbox{}\vspace{.05in} \\
Henk Broer\footnotemark, Richard Cushman\footnotemark, and 
Francesco Fass\`{o}\footnotemark
\end{center}
\addtocounter{footnote}{-2}
\footnotetext{Instituut voor Wiskunde en Informatica, Rijksuniversiteit Groningen,
Blauwborgje 3, 9747 AC Groningen, The Netherlands (email: broer@math.rug.nl)} 
\addtocounter{footnote}{1}
\footnotetext{Faculteit Wiskunde en Informatica, Universiteit Utrecht,
Budapestlaan 6, 3584 CD  Utrecht, The Netherlands (email: cushman@math.uu.nl)}
\addtocounter{footnote}{1} 
\footnotetext{Universit\`a di Padova, Dipartimento di
Matematica Pura e Applicata, Via G. Belzoni 7, 35131 Padova,
Italy (email: fasso@math.unipd.it)}

\begin{abstract}
\noindent
We obtain a global version of the Hamiltonian KAM theorem for invariant 
Lagrangean tori by glueing together local KAM conjugacies 
with help of a partition of unity. In this 
way we find a global Whitney smooth conjugacy between
a nearly-integrable system and an integrable one.
This leads to preservation of geometry, which allows us to 
define all the nontrivial geometric invariants like monodromy or
Chern classes of an integrable system also for near integrable systems.
\end{abstract}

\section{Introduction}

Classical Kolmogorov-Arnold-Moser theory deals with Hamiltonian perturbations 
of an integrable Hamiltonian system and proves the persistence of 
quasi-periodic (Diophantine) invariant Lagrangean tori. In 
\cite{Pos} P\"{o}schel proved the existence of Whitney smooth 
action angle variables on a nowhere dense union of tori having 
positive Liouville measure. See also \cite{Laz,CG}. 
This version of the KAM theorem can be 
formulated as a kind of structural stability restricted to 
this union of quasi-periodic tori. As such it is referred to as 
\emph{quasi-periodic stability} \cite{BHS,BHT}. 
In this context, the conjugacy between the integrable system and its perturbation 
is smooth in the sense of Whitney. 
This paper deals only with the case of Hamiltonian systems with 
Lagrangean invariant tori. \medskip

In this paper we aim to establish a global quasi-periodic stability 
result for fibrations of Lagrangean tori,
by glueing together local conjugacies obtained from the 
classical KAM theorem. This glueing uses a partition of unity \cite{Hir,Spivak}
and the fact that invariant tori of the unperturbed integrable system have 
a natural affine structure \cite{Arnold,Cushman-Bates,Fas}. The global 
conjugacy is obtained as an appropriate convex combination of 
the local conjugacies. This construction is reminiscent of that used to construct 
connections or Riemannian metrics in differential geometry.

\subsection{Motivation}

A motivation for globalizing the KAM theorem is the 
nontriviality of certain torus fibrations in Liouville integrable systems. 
For example, in two degrees of freedom, consider the spherical pendulum. 
Here an obstruction to the triviality of the foliation by 
Liouville tori is given by monodromy, see \cite{Duis,Cushman-Bates}.
A natural question is whether (nontrivial) monodromy also can be defined 
for nonintegrable perturbations of the spherical pendulum. Answering this question is 
of interest in the study of semiclassical versions 
of such classical systems, see \cite{montgomery, cushman-sadovskii}. 
The results of the present paper imply that for an open set of Liouville integrable Hamiltonian systems, 
under a sufficiently small perturbation, the geometry of the fibration is largely preserved by 
a Whitney smooth diffeomorphism. For a geometrical discussion of all the obstructions for 
a toral fibration of an integrable Hamiltonian system to be trivial see \cite{Duis}. \medskip

\noindent \textbf{Remark.} For a similar result in two degrees of freedom 
near focus-focus singularities see \cite{Rink}. We expect that 
generalizations of our results will be valid in the Lie algebra setting of 
\cite{BHS,BHT}.

\subsection{Formulation of the result}

We now give a precise formulation of our result. Consider 
a $2n$-dimensional, symplectic, real analytic manifold 
$(\widetilde{M},\sigma)$ with a surjective real analytic map 
$\widetilde{\pi}: \widetilde{M}\rightarrow \widetilde{B},$
where $\widetilde{B}$ is an $n$--dimensional affine manifold.
The map $\widetilde{\pi }$ is a foliation with 
singularities in the sense of 
Stefan-Sussman \cite{Stefan,Sussmann}. Its  
regular leaves are Lagrangean $n$--tori.
If $B \subseteq  \widetilde{B}$ is the set of regular values of $\widetilde{\pi}$, 
we assume that $M = {\widetilde{\pi}}^{-1} (B)$ is connected. 
Let  $\pi: M\rightarrow B$ be the restriction and corestriction of 
$\widetilde{\pi}$ to ${\widetilde{\pi }}^{-1}(B)$. 
By the Liouville--Arnold integrability theorem \cite{Arnold,Cushman-Bates} 
it follows that for every $b \in B $ there is a neighbourhood
$U^b \subseteq B$ and a symplectic diffeomorphism 
\begin{displaymath}
\varphi:  V^b = \pi^{-1} (U^b) \rightarrow \mT^n \times A^b, \, \, 
m \mapsto (\alpha^b (m), a^b (m)),
\end{displaymath}
with symplectic form $\sum_{j=1}^n da_j^b\wedge d \alpha_j^b$ 
such that $a^b = (a_1^b, a_2^b, \ldots, a_n^b)$ is constant on fibers of $\pi.$ 
As usual we call $(a^b, \alpha^b)$ action angle variables and 
$(V^b, \varphi^b)$ an action angle chart.

Now consider a real analytic Hamiltonian function $H: \widetilde{M} \rightarrow \mR,$ 
which is integral of $\pi $, that is, $H$ is constant on fibers of $\widetilde{\pi}$. 
Then the corresponding Hamiltonian vector 
field $X_H$, defined by $X_{H} \lefthook \sigma = \dee H$, is 
tangent to these fibers. This leads to a push forward vector field
\[
(\varphi^b_\ast X_H) (\alpha,a)  = \sum_{j=1}^n \omega^b_j (a)\, 
                      \frac{\partial}{\partial \alpha_j^b},
\]
on $A^b \subseteq {\R }^n$ with frequency vector 
$\omega^b (a) = (\omega^b_1 (a),\ldots,\omega^b_n (a))$. 
We call $\omega^b : A^b \rightarrow \mR^n$ the local frequency map. \medskip 

We say that $H$ is a \emph{nondegenerate} integral of $\widetilde{\pi},$
if for a collection ${\{ (V^b, \varphi^b) \} }_{b\in B}$ of action angle 
charts whose domains $V^b$ cover $M$, each local frequency map  $\omega^b : A^b \rightarrow \mR^n$ 
is a diffeomorphism onto its image. \medskip 

\noindent
\textbf{Remark.} Using the (local) affine structure on the space of 
actions $B \subseteq \widetilde{B}$, for each $b\in B$ the second derivative of 
$h^b = (H|V^b) \comp ({\varphi }^b)^{-1}$ on $V^b$ is well-defined. By nondegeneracy we just mean 
that $D^2 h^b$ has maximal rank everywhere on 
$V^b$ for each $b \in B$. \medskip

\noindent 
Suppose that $H$ is a nondegenerate integral on $B$. 
If $B'\subseteq B$ is a relatively compact subset of $B$, 
that is, the closure of $B'$ is compact, then there is a finite 
subcover ${\{ U^b \} }_{b\in \mathcal{F}}$ of $\{ U^b \}_{b\in B}$ such that 
for every $b \in \mathcal{F}$ the frequency map $\omega^b $ is a diffeomorphism onto its image. 
Accordingly, let $M' = \pi^{-1} (B')$ and consider the corresponding bundle
$\pi': M' \rightarrow B'.$ 
We now state our main result. \medskip 

\begin{description}
\item[Theorem 1.]{\sc (KAM-global)} \label{th:globalKAM}
Let $(M,\sigma)$ and $H: M \rightarrow \mR$ be as above 
with $B'\subseteq B$ open and relatively compact. 
Let $F: M \rightarrow \mR$ be a real analytic function.\footnote{Generally 
the perturbation $H+F$ is not integrable.}
If $F|_{M'}$ is sufficiently 
small in the compact-open topology 
on holomorphic extensions, then there exists a subset $C \subseteq B'$ and 
a map $\Phi: M \rightarrow M$ with the following properties. 
\begin{enumerate}
\item The subset $C$ is a nowhere dense set of large measure in $B'$; 
\item The map $\Phi$ is a $C^\infty$ diffeomorphism onto its image and is
close to the identity map in the $C^\infty$ topology;
\item The restriction $\widehat{\Phi} = \Phi|_{\pi^{-1}(C)}$ 
conjugates $X_H$ to $X_{H+F}$, that is, 
\begin{displaymath}
\widehat{\Phi}_\ast X_H = X_{H + F}.
\end{displaymath} 
\end{enumerate}
\end{description}

\clearpage
\noindent
\textbf{Remarks.} 
\begin{enumerate}
\item The map $\Phi$ generally is not symplectic. Since its derivatives are bounded, 
the bundle $\pi': M'\rightarrow B'$ and its perturbation are diffeomorphic. 
Moreover, the affine structure on the Diophantine tori is preserved.
\item  The measure of $\Phi (\pi^{-1}(C))$, which is the union  
of perturbed tori, is large.
\item The restriction to real-analytic systems with the 
compact-open topology on holomorphic extensions is not essential.
Indeed, there are \linebreak 
direct generalizations of theorem 1 to the world of $C^k$-systems 
endowed with the (weak) Whitney  $C^k$-topology  
for $k$ sufficiently large, see \cite{Hir,BTan}.
For similar $C^k$-versions of the local KAM theorem, see \cite{Pos,BHT,BHS}.
\end{enumerate}

\section{Proof}
This section is devoted to proving theorem 1.
The idea is to glue together the local conjugacies (obtained
by the standard local KAM theorem \cite{Pos,BHT}) with 
help of a partition of unity. 

\subsection{The local KAM theorem}
We need a precise formulation of the local KAM theorem.
The phase space here is $\mT^n \times A,$ where 
$A\subseteq \mR^n$ is bounded, open, connected
and $(\alpha, a) = (\alpha_1,\ldots\alpha_n,a_1,\ldots,a_n)$
is a set of action angle coordinates with 
symplectic form $\sum_{j = 1}^n \dee a_j \wedge \dee \alpha_j.$
Suppose that we are given a real analytic Hamiltonian function 
$h: \mT^n \times A \rightarrow \mR,$  which is integrable, that is,  
$h$ does not depend on the angle variables $\alpha$. 
We shall now formulate a few concepts we need in the proof of theorem 
$1$.  

\medskip
The Hamiltonian vector 
field $X_h$ corresponding to $h$ has the form
\begin{displaymath}
X_h (\alpha, a) = 
\sum_{j =1}^n \omega_j (a) \frac{\partial}{\partial \alpha_j}, 
\end{displaymath}
with $\omega (a) = \partial h / \partial a.$ The nondegeneracy 
assumption now means that the corresponding frequency map 
$\omega: a\in A \mapsto \omega (a) \in \mR^n$ is a 
diffeomorphism onto its image.

\medskip
Also we need the concept of Diophantine frequencies.
Let $\tau > n-1$ be a fixed constant and let $\gamma > 0$ be a `parameter'. 
Let 
\[
{\rm D}_{\gamma} (\mR^n) = 
\{\omega \in \mR^n \mid 
|\langle \omega,k \rangle| \ge \gamma |k|^{-\tau}, 
\mbox{ for all } k \in \mZ^n \setminus \{0\} \}.
\]
Elements of ${\rm D}_{\gamma}(\mR^n) $ are called \emph{Diophantine} frequency vectors.
Consider the set $\Gamma = \omega (A).$ For $\tilde{\gamma} > 0$ also consider the shrunken version 
$\Gamma_{\tilde{\gamma}} \subseteq \Gamma$ given by
\[
\Gamma_{\tilde{\gamma}} = \{\omega \in \Gamma = \omega (A) \mid 
\mathrm{dist} (\omega, \partial \Gamma) > \tilde{\gamma}\} .
\]
Let ${\rm D}_\gamma (\Gamma_{\tilde{\gamma}}) = \Gamma_{\tilde{\gamma}} \cap {\rm D}_\gamma (\mR^n).$
From now we take $\tilde{\gamma} = \gamma$ sufficiently small to ensure that 
${\rm D}_\gamma (\Gamma_{\gamma})$  is a nowhere dense set of positive measure,
 compare \cite{Pos,BHS,BHT}. Recall that this measure tends to full measure 
in $\Gamma$ as $\gamma \downarrow 0.$
Finally define the shrunken domain $A_\gamma = \omega^{-1} (\Gamma_\gamma) \subseteq A$ 
as well as its nowhere dense counterpart 
${\rm D}_\gamma (A_\gamma) = \omega^{-1} ({\rm D}_\gamma(\Gamma_\gamma)) \subseteq A_\gamma ,$
which again has large measure in $A$ for $\gamma$ small.
 \medskip

We now perturb $h$ to a real analytic Hamiltonian function 
\begin{displaymath}
h + f: \mT^n \times A \rightarrow \mR,
\end{displaymath} 
which we assume is holomorphic on the 
complexified domain
\[
D_{\varrho,\kappa} = (\mT^n + \kappa) \times (A + \varrho),
\]
for positive constants $\kappa$ and $\varrho < 1$. Here 
\[
A+\varrho = \{y\in \mC \mid |y-a| \le \varrho \mbox{ for some } a \in A \}.
\]
Similarly we define $\mT^n + \kappa \subseteq \mC /(2\pi \mZ).$ \medskip 

We now are ready to formulate the standard local KAM theorem.

\begin{description}
\item[Theorem 2.] {\sc (KAM-local)} \mbox{\hspace{4pt}\cite{Pos,BHT}} Suppose that the integrable real 
\linebreak 
analytic Hamiltonian function $h: \mT^n \times A\rightarrow \mR$ is nondegenerate and that both $h$ and 
its real analytic 
perturbation $h + f$ have a holomorphic extension 
to the domain $D_{\varrho,\kappa} = (\mT^n + \kappa) \times (A + \varrho),$
with positive constants $\kappa$ and $\varrho < 1.$ 
Then 
\begin{itemize}
\item[1)] there exists a positive constant $\delta$, \emph{independent} 
of $\gamma$, $\varrho $, and of the domain $A$; 
\item[2)] if
$|f|_{D_{\varrho,\kappa }} \le {\gamma}^2 \delta $, where $|\, \, |_{D_{\varrho,\kappa } }$ 
is the supremum-norm on $D_{\varrho ,\kappa }$, 
then there is a map $\Phi^A: \mT^n \times A \rightarrow \mT^n \times A $ 
with the following properties 
\begin{itemize}
\item[a)] $\Phi^A$ is a $C^\infty$ diffeomorphism onto its image, 
      which is real analytic in the angle variables $\alpha$ and is 
      close to the identity map in the $C^\infty$ topology;
\item[b)] The map $\widehat{\Phi^A} = \Phi^A|_{\mT^n \times {\rm D}_\gamma (A_\gamma)}$ 
      conjugates $X_h$ to $X_{h+f}$, that is, ${\widehat{\Phi^A}}_\ast X_h 
= X_{h + f}$.
\end{itemize}
\end{itemize}
\end{description}
\noindent \textbf{Remarks.} 
\begin{enumerate}
\item Observe that $\Phi^A $ maps $A_{\gamma }$ into $A$. 
\item  In general, the map $\Phi^A$ is not symplectic. Notice 
      that theorem 2 is stated in its stability form \cite{BHT,BHS}, namely, 
      the theorem asserts that the integrable system is     
      \emph{quasi-periodically stable}.
\item The nowhere dense set ${\rm D}_\gamma (\mR^n)$ has a raylike structure. 
      Whenever $\omega \in {\rm D}_\gamma (\mR^n)$, we also have   
      $s \omega \in {\rm D}_\gamma (\mR^n)$ for $s \ge 1$. In the ray direction 
      of $A_\gamma$ the map $\Phi^A$ is real analytic. 
      Moreover, analytic dependence on extra parameters is preserved 
      in the following sense: if this analytic dependence holds for $h+f$, 
      then it also holds for $\Phi^A$.
\item In applications, given the perturbation 
      the parameter $\gamma$ is chosen as small as possible. 
      In particular, if the perturbation is of the form  $h+\varepsilon f$, 
      then we can take $\gamma = \mathrm{O}(\sqrt{\varepsilon})$.    
\end{enumerate}

\subsection{Application of the local KAM theorem}

Let us return to the global problem of perturbing $H$ to $H+F$ on $M' = \pi^{-1}(B').$ 
Recall that we have an atlas $\{(V^b,\varphi^b)\}_{b\in B}$ of action angle charts on $M,$ where
$V^b = \pi^{-1} (U^b) = ({\varphi^b})^{-1} (\mT^n \times A^b).$ For each $b\in B$ there
exists a constant $\gamma^b > 0$ such that the shrunken domain $V_{\gamma^b}^b \subseteq V^b$
defined by $V_{\gamma^b}^b = ({\varphi^b})^{-1} (\mT^n \times A^b_{\gamma^b})$
is open and nonempty. This implies that $\{V^b_{\gamma^b}\}_{b\in B}$ still is a covering of $M.$
Given any choice of an open relatively compact set $B' \subseteq B,$ 
there exists a finite subcover ${\{ V_{\gamma^j}^j \} }_{j\in {\cal J}}$ 
of $M' = \pi^{-1} (B').$  Note that $V_{\gamma^j}^j \subseteq V^j,$ $j \in {\cal J}.$
In each action angle chart $(V^j, {\varphi }^j)$ where 
$\varphi^j : V^j = \pi^{-1} (U^j) \rightarrow \mT^n \times A^j$, we have 
a local perturbation problem. Indeed, defining 
$h^j = H \comp (\varphi^j)^{-1}$ and  $f^j = F \comp (\varphi^j)^{-1}$, 
we are in the setting of the local KAM theorem 2 on the phase space
$\mT^n \times A^j$. Accordingly we are given positive constants $\kappa = \kappa^j$ and 
$\varrho = \varrho^j$ and we specify $\gamma = \gamma^j$ as just obtained. \medskip

In this setting the nowhere dense set $C \subseteq B'$ is defined by pulling back 
the sets ${\rm D}_{\gamma^j}(A^j_{\gamma^j}) \subseteq A^j_{\gamma^j}\subseteq A^j$ 
using the second component of the chart map $\varphi^j$. On overlaps we just take the intersection 
of finitely many domains $V^j_{\gamma^j} \subseteq V^j.$ Note that for small $\gamma^j,$ $j\in{\cal J},$
the nowhere dense set $\pi^{-1}(C) \subseteq M'$ has large measure. 
(Compare with the first claim of the global KAM theorem 1).

Therefore, given the appropriate (equicontinuous) smallness conditions on 
the collection of functions ${\{ f^j \} }_{j\in {\cal J}}$, by the local KAM theorem 2 
we obtain diffeomorphisms 
$\Phi^j = \Phi^{A^j}: \mT^n \times A^j \rightarrow \mT^n \times A^j,$ such that 
restricted to ${\rm D}_{\gamma^j} (A^j_{\gamma^j}),$ we have 
\[
(\widehat{\Phi^j})_\ast X_{h^j} = X_{h^j + f^j}.
\]
 
From the above construction we deduce that each of the Diophantine invariant tori 
$T$ of $X_H$ under consideration is diffeomorphic by $(\varphi^j)^{-1}\comp \Phi^j \comp \varphi^j$ 
to an invariant torus $T'$ of $X_{H + F}.$ This map 
conjugates the quasi-periodic dynamics of $X_H$ and $X_{H+F},$
where the frequency vector is invariant. Thus, by 
nondegeneracy of $H$, the correspondence which associates the 
torus $T$ to the torus $T'$ is unique, and hence is 
independent of the index $j\in {\cal J}$ used to define the conjugation. 
In other words,  on the union  $\pi^{-1}(C)$ of Diophantine tori
the action coordinates match under the transition maps
$\varphi^j \comp (\varphi^i)^{-1}.$

\medskip  
At the level of the angle coordinates the matter of matching is quite different. 
First, in cases where the $n$-torus bundle $\pi: M\rightarrow B$ is nontrivial,
generally the angle coordinates on $\pi^{-1}(C)$ do not have to 
match under the transition maps $\varphi^j \comp (\varphi^i)^{-1}.$
Second, the angle components of local KAM conjugacies $\Phi^i$ and $\Phi^j$
generally do not have to match either.

\subsection{Partition of Unity}

To overcome the two matching problems at the level of angles as mentioned above, 
we introduce a partition of unity.
\begin{description}  
\item[Lemma 1. ({\sc partition of unity}) ]\label{th:PoU}
Subordinate to the covering  ${\{ V^b_{\gamma^b} \} }_{b\in B'}$ of $M'$ 
there is a partition of unity  
$\{ (V^j_{\gamma^j}, \xi^j) \}_{j \in {\cal J}}$ by $C^\infty $ functions  
$\xi^j : M \rightarrow \mR$ such that, for every $j \in {\cal J}$ we have  
\begin{enumerate}
\item the support of $\xi^j$ is a compact subset of $V^j_{\gamma^j}$;
\item the function $\xi^j$ takes values in the interval $[0,1];$ 
\item the function $\xi^j$ is constant on the fibers of the 
bundle projection $\pi': M' \rightarrow B'$
\item $\sum^N_{j=1} \xi^j \equiv 1$ on $M'$.  
\end{enumerate}
\end{description}

\noindent
The lemma follows by carrying out the standard partition of unity construction 
\cite{Hir,Spivak} on $B'$ and pulling everything back by $\pi.$

\subsection{Glueing}

For each $j \in {\cal J}$ consider the map
\[
\Psi^j = (\varphi^j)^{-1}\comp (\Phi^j)^{-1} \comp \varphi^j,
\]
which takes a perturbed torus $T'$ to its unperturbed (integrable) counterpart $T$ and 
thereby conjugates the quasi-periodic dynamics. 
Note that the local KAM conjugagy $\Phi^j$ maps $\mT^n\times A^j_{\gamma^j}$ into 
$\mT^n\times A^j.$ Therefore $(\Psi^j)^{-1}$ maps $V^j_{\gamma^j}$ into $V^j.$

We now use the fact that the integrable tori $T$ have a natural affine structure,
see the proof of the Liouville-Arnold theorem \cite{Arnold, Cushman-Bates}.
Indeed, this structure is induced by the transitive $\mR^n$-action
\[
     G : \mR^n \times T \rightarrow T, \ \ 
     ((t_1,\ldots,t_n),m) \mapsto ({g }_1^{t_1}\comp \cdots \comp  
{g }_n^{t_n} (m)),
\]
where  ${g }_\ell^t,$  for $\ell =1,\ldots ,n,$ is the flow of the Hamiltonian vector field 
associated to the action $a_\ell^j,$ $j= 1,\ldots , N$.
We make the following assertions about the chart overlap maps 
on $M'.$\footnote{Note that for the tori $T$ and $T'$
the affine structure also is determined by the quasi-periodicity of the flow.} 

\begin{description}
\item[Lemma 2. ({\sc overlap})]\label{th:overlap}
For a perturbation function $F,$ sufficiently close to $0$ in the compact-open topology on 
holomorphic extensions,  on the overlaps $V^i\cap V^j$ 
\begin{enumerate}
\item The transition map $\Psi^i \comp (\Psi^j)^{-1}$
       is close to the identity in the $C^\infty$ topology;
\item For each integrable $n$-torus $T \subseteq V^i_{\gamma^i} \cap V^j_{\gamma^j},$
       which by $(\Psi^i)^{-1}$ and  $(\Psi^j)^{-1}$
       is diffeomorphic to $T',$
       the maps $\Psi^i|_{T'}$ and $\Psi^j|_{T'}$ only differ by a small
       translation. 
\end{enumerate}
\end{description}
\noindent \textbf{Proof.}
The first item directly follows from the equicontinuity condition regarding
the size of the ${\{ f^j \} }_{j\in {\cal J}}$. The second item follows 
from the former, see \cite{Arnold,Cushman-Bates,Fas}. \hfill$\square$ \medskip 

\noindent \textbf{Remark.} On an overlap  $T \subseteq V^i_{\gamma^i}\cap V^j_{\gamma^j},$ 
$i\ne j,$ consider the angle  components  $\alpha^i$ and  $\alpha^j$ of the 
chart maps $\varphi^i$  and $\varphi^j,$ both taking values
in the standard torus $\mT^n.$ Here $\alpha^j  = S_{i,j} \, \alpha^i + c_{i,j}$ 
for $S_{i,j}\in Sl(n,\mZ)$ and $c_{i,j}\in \mT^n$, see \cite{Arnold,Cushman-Bates,Fas}. 
In the case of nontrivial monodromy we do not always have $S_{i,j} = {\rm Id},$
since the transition maps $\varphi^i \circ (\varphi^j)^{-1}$ are not all close 
to the identity map. We have avoided this problem by considering the maps 
$\Psi^i$ instead. \medskip 

\noindent \textbf{Proof of theorem 1.}
The global conjugacy of theorem 1 is obtained 
by taking the sum
\begin{equation}
\Phi^{-1} = \sum_{j=1}^N \xi^j  \Psi^j,
\label{eq-star} 
\end{equation}
using the partition of unity of lemma 1.\medskip

Equation (\ref{eq-star}) has to be interpreted {\sc fibrewise} as follows.
In each integrable torus $T,$ by lemma 2, the maps 
$ \Psi^i$ and $ \Psi^j$ only differ by a small translation.
The finite convex combination is then a conjugacy, which is globally well-defined. 
The fact that for small $F$ the map $\Phi$ is $C^\infty$-close to the identity map, 
follows by Leibniz' rule.
This proves the global KAM theorem 1. 
\hfill$\square$

\begin{figure}[htp]
\begin{picture}(140,45)
\put(13,-50){\resizebox{10 cm}{!}{\includegraphics{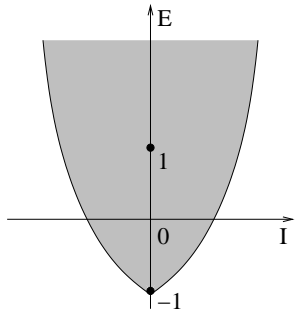}}}
\end{picture}
\caption{Range of the energy--momentum map of the spherical
pendulum. \label{fig:energymomentum}}
\end{figure}

\section{An application}

Applying theorem 1 largely amounts
to checking the global nondegeneracy of the Hamiltonian function $H.$ \medskip

As an example let us consider the spherical pendulum \cite{Duis,Cushman-Bates}. 
The configuration space is the 2-sphere
$S^2 = \{q\in \mR^3 \mid \langle q,q\rangle = 1 \}$ 
and the phase space its cotangent bundle
$\widetilde{M} = 
T^{\ast }S^2 =  \{ (q,p)\in \mR^6 \mid \langle q,q \rangle = 1 \, \, 
\mathrm{and} \, \,  \langle q,p\rangle = 0\}.$
The energy-momentum map ${\cal EM}: T^{\ast }S^2 \rightarrow \mR^2$, given by 
\[
{\cal EM}\, (q,p) = (I,E) = 
\left(q_1p_2 - q_2 p_1, H(q,p) = \onehalf \langle p,p \rangle + q_3 \right),
\]
leads to a Lagrangean fibration, whose generic fiber is a 
$2$-torus. The image $\widetilde{B}$ of ${\cal EM}$
is the closed part of the plane lying in between the two curves
meeting at a corner, see figure 1. 
The set of singular values of $\mathcal{EM}$ consists of the 
two boundary curves and the points $(I,E) = (0,\pm 1).$ The latter correspond
to the equilibria  $(q,p) = (0,0,\pm 1,0,0,0).$ The boundary
curves correspond to the horizontal periodic motions of the pendulum.
The set $B$ of regular ${\cal EM}$-values therefore consists of the 
interior of $\widetilde{B}$ minus  the point $(I,E) = (0,1),$ corresponding 
to the unstable equilibrium point $(0,0,1,0,0,0)$. The point 
$(I,E) = (0,1)$ is  the centre 
of the nontrivial monodromy. The corresponding fiber 
${\mathcal{EM}}^{-1}(1,0)$ is a once pinched $2$-torus.

On $B$ one of the two frequency maps is single valued and the other
multi-valued \cite{Duis,Cushman-Bates}. Horozov \cite{horozov} 
established global nondegeneracy of $H$ on $B.$ This implies that 
the global KAM theorem 1 can be applied to any relatively compact open subset 
$B' \subseteq B.$ From this we deduce that the integrable dynamics on the 
$2$-torus bundle ${\mathcal{EM}}': M' \rightarrow B'$ of the spherical 
pendulum is quasi-periodically stable. 
This means that any sufficiently small perturbation of the spherical pendulum, integrable or not, 
has an equivalent KAM $2$-torus bundle which has the same geometry as 
the bundle ${\mathcal{EM}}'$. In particular,
the torus bundle can be said to have nontrivial monodromy. 

Asymptotics near the horizontal periodic solutions and both the 
equilibria lead to the following considerations. (Compare this with 
the problem of small twist \cite{Pos}.) 
The parameter $\gamma $ as it occurs in the Diophantine conditions,
can be taken small in an appropriate way, see \cite{BHT,BHS},
thereby creating Lebesgue density points of quasi-periodicity
at the horizontal periodic solutions and at both equilibria. \medskip 

\noindent \textbf{Remarks.} 
\begin{enumerate}
\item From the asymptotic considerations of Rink \cite{Rink} near focus-focus singularities 
     in Liouville integrable Hamiltonian systems of $2$ degrees of freedom, 
     the above results regarding the spherical pendulum 
     follow in the case where $B'$ is a small annular region around  $(0,1).$ 
     This is sufficient for the existence of nontrivial monodromy in nearly integrable perturbations of the 
     spherical pendulum. 
\item The present global approach does not depend on knowing the 
      precise geometry of the integrable approximation. It clearly 
      applies near focus-focus singularities.
\item A straightforward generalization of the global KAM theorem 1
     exists in all cases of (local) quasi-periodic stability considered in 
     \cite{BHT}. In \cite{BHT}, partly following \cite{Moser}, 
     a general Lie algebraic approach is taken, 
     which leads to various versions of the (local) KAM theorem. For instance, 
     the general (disspative) case, the volume preserving case, etc, see 
     \cite{BHS}. This includes the case of lower dimensional isotropic tori 
     in Hamiltonian systems. Often external parameters are needed for 
     quasi periodic stability. In this case the torus bundles live in the 
     product of phase space and parameter space. For a reversible analogue 
     see \cite{BH}. We expect that the present glueing construction directly 
     carries over to global torus bundles in these settings.
\item The present paper opens the possibility to define geometric concepts like
    monodromy, Chern classes, etc.\ (compare \cite{Duis}), for KAM torus bundles in 
    nearly integrable systems, by integrable approximation. It is a
    quite different matter how to define such concepts in a more direct way,
    say by  construction of certain appropriate Whitney smooth sections.

\end{enumerate}

\section*{Acknowledgements}
The authors are grateful to Bob Rink, Michael Sevryuk, Carles  Sim\'o, 
and to Floris Takens for helpful discussions. 
Also the first author (HWB) acknowledges the hospitality of the 
Universitat de Barcelona. The authors were partially supported by 
European Community funding for the Research and Training Network 
MASIE (HPRN-CT-2000-00113).

\end{document}